\documentclass[1p]{elsarticle}

\usepackage{amsfonts}
\usepackage{amssymb}
\usepackage{bbm}
\usepackage{amsmath}
\DeclareMathOperator{\Tr}{Tr}
\DeclareMathOperator{\Hess}{Hess}

\usepackage{graphicx}
\usepackage{array}
\usepackage{comment}

\usepackage{multirow}
\usepackage{float}
\usepackage{algorithm2e}
\RestyleAlgo{ruled}

\usepackage{xcolor}
\usepackage[colorlinks=true, allcolors=blue]{hyperref}

\journal{Computers and Mathematics with Applications}

\begin{document}

\begin{frontmatter}

\title{A Deep-Genetic Algorithm (Deep-GA) Approach for High-Dimensional Nonlinear Parabolic Partial Differential Equations}

\author[address1]{Endah RM Putri\corref{mycorrespondingauthor}}
\cortext[mycorrespondingauthor]{Corresponding author}
\ead{endahrmp@matematika.its.ac.id}

\author[address2,address1]{Muhammad L Shahab}
\author[address1]{Mohammad Iqbal}
\author[address1]{Imam Mukhlash}
\author[address1]{Amirul Hakam}
\author[address3]{Lutfi Mardianto}
\author[address2]{Hadi Susanto}

\address[address1]{Department of Mathematics, Faculty of Science and Data Analytics, Institut Teknologi Sepuluh Nopember, Jl. Raya ITS, Sukolilo, Surabaya, 60111, Indonesia}
\address[address2]{Department of Mathematics, Khalifa University, Abu Dhabi, PO Box 127788, United Arab Emirates}
\address[address3]{Department of Mathematics, Institut Teknologi Sumatera, Jl. Terusan Ryacudu, Way Hui, Jati Agung, Lampung Selatan, 35365, Indonesia}

\begin{abstract}
We propose a new method, called a deep-genetic algorithm (deep-GA), to accelerate the performance of the so-called deep-BSDE method, which is a deep learning algorithm to solve high dimensional partial differential equations through their corresponding backward stochastic differential equations (BSDEs). Recognizing the sensitivity of the solver to the initial guess selection, we embed a genetic algorithm (GA) into the solver to optimize the selection. We aim to achieve faster convergence for the nonlinear PDEs on a broader interval than deep-BSDE. Our proposed method is applied to two nonlinear parabolic PDEs, i.e., the Black-Scholes (BS)
equation with default risk and the Hamilton-Jacobi-Bellman (HJB) equation. We compare the results of our method with those of the deep-BSDE and show that our method provides comparable accuracy with significantly improved computational
efficiency.
\end{abstract}

\begin{keyword}
high dimensionality \sep non-linear equations \sep genetic algorithm\sep backward stochastic differential equation
\end{keyword}

\end{frontmatter}

\section{Introduction}

Many natural phenomena and complex systems are modeled by parabolic partial
differential equations (PDEs), particularly in financial industry such
as derivative pricing and portfolio optimization models. Nonlinear parabolic
PDEs such as the Black-Scholes (BS) with default risk and the Hamilton-Jacobi-Bellman
(HJB) equations are popular examples of PDEs in the field
\cite{han2018solving}, while other examples of nonlinear parabolic PDEs can be found
in
\cite{bernhart2012swing,esmaeeli2018american,sun2021bsde,cordoni2014backward}.
Nonlinearity in the models aims to incorporate phenomena like default
risks \cite{liao2009high}, transaction costs \cite{kohler2010pricing},
and stochastic volatilities \cite{goudenege2020machine}. As most nonlinear
models have no analytical solution, approximation or numerical methods
are required. Adding high-dimensionality to the nonlinear financial models
increases the complexity of the solution process
\cite{weinan2021algorithms}. Moreover, high-dimensional nonlinear models
bring the so-called curse of dimensionality (CoD) as extensive and complex
computations grow exponentially. Classical grid based methods
i.e., finite difference \cite{kim2020finite,miyamoto2021pricing}, finite
element \cite{zhang2015efficient}, or finite volume methods
\cite{moroney2006finite}, fail in high dimensional cases and work only in the lower one or two dimensions \cite{weinan2021algorithms}.

Monte Carlo method is a
powerful alternative to obtain accurate solutions of high dimensional
problems where CoD exists. In general, the method employs random simulations
to generate all possible solutions and determines the optimal one. It is then applied to linear PDEs in the form of Kolmogorov backward problems. After a slow development of the algorithm, a multilevel Picard approximation
method (MLP) made a break through \cite{weinan2021algorithms} as a nonlinear Monte Carlo approximation method.
Under certain assumptions, the MLP was proven analytically to solve
CoD in nonlinear PDEs \cite{hutzenthaler2019multilevel}.%

The so-called deep-BSDE algorithm has been recognized as the first deep learning
based method \cite{han2018solving,han2017deep,weinan2017deep} in solving
high dimensional semi-linear parabolic PDEs. The method represents the PDEs, through their corresponding backward stochastic differential equations (BSDEs), as a deep learning (stochastic optimization) problem that finds an optimal initial guess satisfying a given terminal condition.
The deep learning based method is applicable to  various high-dimensinal complex problems, such as an optimal stopping problem
\cite{becker2019deep,becker2021solving}, a fixed point problem
\cite{chan2019machine}, a many-electron Schr\"{o}dinger equation as an eigenvalue problem \cite{han2020solving}, etc. Accordingly, the method has been extended in several directions, such as into least square based deep learning  method for not only semi-linear PDEs
but also general PDEs \cite{SIRIGNANO20181339}, deep-BSDE using second
order BSDEs for fully nonlinear PDEs \cite{beck2019machine}, deep splitting
method for parabolic PDEs which separate the linear and nonlinear terms
to save the computational time \cite{beck2021deep}, etc.

The aforementioned studies commonly employ a predetermined initial guess, which we assume comes from a trial and error process. However,
we argue that appropriately selecting the initial guess is crucial and can contribute to a faster convergence, as we will demonstrate in Section~\ref{sec4} below. Here, we propose to employ a genetic algorithm (GA) to achieve the advantages. Pan et al.\ \cite{pan2014efficient} found that
employing both GA and gradient descent algorithms
enables a faster global search for optimal weights in neural networks.
They also conclude that the use of GA outperforms
other evolutionary algorithms, such as a differential evolution and a particle
swarm optimization. Ding et al.\ \cite{ding2011optimizing} showed that combining
GA and a gradient descent algorithm produces better
solutions and has a good and stable performance. As an example to applications
in PDEs, a hybrid deep learning and GA for data-driven
discovery of PDEs is presented in \cite{xu2020dlga}. Moreover, GA is applicable to many areas
\cite{bouktif2018optimal,kilicarslan2021hybrid,kwon2021evolutionary,kalsi2018dna,balaha2022hybrid,skandha2022novel}.

Instead of relying on a fixed initial guess, we present a deep-genetic
algorithm (deep-GA) method that fuses GA with the Adam
optimizer to address the original deep learning problem discussed in
\cite{han2018solving}. The learning process of the deep-GA, involving weight
updates, is divided alternately between GA and the Adam
optimizer. The primary objective of this research is to develop
an efficient method for solving the deep learning problem emerging
from high-dimensional nonlinear parabolic PDEs.

This paper is organized as follows. Section~\ref{sec2} presents the governing
equations including the set-up of the PDEs and BSDEs of the nonlinear
BS and HJB equation. We will also revisit
the deep-BSDE method. Section~\ref{sec3} discusses the set-up of GA in the BSDE equation and the method in
\cite{han2018solving}. Simulation results and discussions are provided
in Section~\ref{sec4}. Finally, conclusion is given in Section~\ref{sec5}.

\section{Mathematical equations}
\label{sec2}

Solving high-dimensional nonlinear parabolic PDEs is known for their CoD. We consider a general form of the class of equations and its relation to BSDE \cite{han2018solving}. We will give a brief description about the solution process of the typical PDE.

In this paper, we consider two PDEs, namely a
nonlinear BS equation with default risk, and the HJB equation
that arises when considering a classical linear-quadratic Gaussian (LQG)
control problem \cite{han2018solving}. These two equations will serve as
test beds for our method. For the reader's convenience, a concise explanation
of the two equations is provided in this section.

\subsection{Backward Stochastic Differential Equations (BSDEs)}

We are interested the following class of nonlinear parabolic PDEs
\cite{han2018solving}
\begin{align}
\label{pde1}
\nonumber
\frac{\partial u}{\partial t}(t,x)+\frac{1}{2}\Tr \left (\sigma
\sigma ^{\text{T}}(t,x)(\Hess _{x} u)(t,x)\right )&+\nabla u(t,x)
\cdot \mu (t,x)
\\
&+ f\left (t,x,u(t,x),\sigma ^{\text{T}} \nabla u(t,x)\right )=0, 
\end{align}
with a terminal condition $u(T,x)=g(x)$ for $x=\{x_{1},x_{2}, \cdots , x_{d}\}$ and time $t$. A volatility function $\sigma (t,x)$ is a $d\times d$ matrix-valued function, $\mu (t,x)$ is a known vector-valued function, $\sigma ^{\text{T}}(t,x)$ is the transpose of $\sigma (t,x)$, $\Hess _{x} u$ denotes the Hessian of function $u(t,x)$ with respect to $x$, $\Tr$ is the matrix trace, and
$f(\cdot )$ is a given nonlinear function. The aim is to find the
solution $u(0,X_{0})$ at $t=0$ and
$x=X_{0} \in \mathbbm{R}^{d}$.

Let $\{ X_{t} \}_{t\in [0,T]}$ be a stochastic process in $d$ dimension
which satisfies
\begin{equation}
X_{t}=X_{0}+\int _{0}^{t}\mu (s,X_{s})ds+\int _{0}^{t}\sigma (s,X_{s})dW_{s},
\label{eq2}
\end{equation}
where $\{ W_{t} \}_{t\in [0,T]}$ is a Brownian motion (Wiener process).
The solution of Eq. (\ref{pde1}) satisfies the following BSDE
\cite{pardoux1992backward}
\begin{align}
\label{3}
\nonumber
u(t,X_{t})=&\;u(0,X_{0})-\int _{0}^{t}f\left (s,X_{s},u(s,X_{s}),
\sigma ^{\text{T}}(s,X_{s})\nabla u(s,X_{s})\right )ds
\\
&+\int _{0}^{t}[\nabla u(s,X_{s})]^{\text{T}}\sigma (s,X_{s})dW_{s}.
\end{align}

To solve Eq. \eqref{3} for $u(0,X_{0})$ numerically, one normally employs
a temporal discretization to partition the time interval $[0,T]$ into
$0=t_{0}<t_{1}<...<t_{N}=1$. A simple Euler scheme is then applied to obtain
\begin{equation}
\label{4}
X_{t_{n+1}}\approx X_{t_{n}} + \mu (t_{n},X_{t_{n}}) \Delta t_{n}+
\sigma (t_{n},X_{t_{n}})\Delta W_{n}
\end{equation}
and
\begin{align}
\label{5}
\nonumber
u(t_{n+1},X_{t_{n+1}}) \approx & \; u(t_{n},X_{t_{n}}) -f(t_{n},X_{t_{n}},u(t_{n},X_{t_{n}}),
\sigma ^{\text{T}}(t_{n},X_{t_{n}})\nabla u(t_{n},X_{t_{n}}))\Delta t_{n}
\\
&+[\nabla u(t_{n},X_{t_{n}})]^{\text{T}}\sigma (t_{n},X_{t_{n}})
\Delta W_{n}
\end{align}
where
\begin{equation}
\label{6}
\Delta t_{n}=t_{n+1}-t_{n}, \;\;\;\;\; \Delta W_{n}=W_{t_{n+1}}-W_{t_{n}}
\end{equation}
for $n=1,\dots ,N-1$. 

\subsection{Non-linear Black-Scholes equation with default risk}\label{NLBS}

A standard BS equation is a linear parabolic PDE that becomes nonlinear when default risk is incorporated. The nonlinearity in these PDEs leads to challenges in obtaining analytical
solutions. The complexity of these PDEs increases with the presence of
multi-assets in financial products combined with the nonlinear nature
of the problem. The BS equation for a multi-asset European contingent
claim with default risks is explained in the following.

Let $u(t,x)$ be the fair price function of a multi-asset European contingent
claim with default risks in $d$ dimension and time $t$, where the
$d$-dimensional underlying assets are represented by
$x=\{x_{1},x_{2}, \dots , x_{d}\}$. We refer to a general form of a nonlinear
PDE in Eq. (\ref{pde1}) for the financial contract
model. We define $\sigma (t,x)=\hat{\sigma}\,\text{diag}(x)$ and
$\mu (t,x)=\hat{\mu}x$ where $\hat{\sigma}$ and $\hat{\mu}$ are known constant
volatility and drift of the contingent claim return, respectively.

The nonlinear part of the equation is represented by a default risk function
$Q(u(t,x))$ which relies on the fair price of a contingent claim function
$u(t,x)$. The terms $v^{h}$ and $v^{l}$ denote rate of default risks for
high and low thresholds of the contingent claim's fair price
$u(t,x)$, respectively. It is assumed that
$v^{h}<v^{l}, \gamma ^{h}>\gamma ^{l}$, and a recovery rate
$\delta \in [0,1)$. The function is modeled as a first jump of the Poisson
process in a possible default within three regions: high risk
$(-\infty , v^{h})$, moderate risk $[v^{h}, v^{l}]$, and low risk
$(v^{l},\infty )$. The value of default risk $Q(u(t,x))$ is equal to
$\gamma ^{h}$ if the contingent claim price $u(t,x)<v^{h}$ or falls in
the high risk region. Conversely, $u(t,x)>v^{l}$ or falls in the low risk
region, the value of $Q(u(t,x))$ is $\gamma ^{l}$. In region
$[v^{h}, v^{l}]$, the value of $Q(u(t,x))$ is classified as moderate risk
and is obtained by an interpolation of the default risk $Q(u(t,x))$ in
the high risk and low risk regions.

The nonlinear BS equation with default risk as the model of
the contingent claim can be written as
\begin{align}
\nonumber
\frac{\partial u}{\partial t}(t,x)+\hat{\mu}x\cdot \nabla u(t,x)&+
\frac{\hat{\sigma}^{2}}{2}\sum _{i=1}^{d}|x_{i}|^{2}
\frac{\partial ^{2}u}{\partial x_{i}^{2}}(t,x)
\\
&+f\left (t,x,u(t,x),\sigma ^{\text{T}} \nabla u(t,x)\right )=0
\label{Eq.2}
\end{align}
where
\begin{align*}
f\left (t,x,u(t,x),\sigma ^{\text{T}} \nabla u(t,x)\right )=-(1-\delta )Q
\left (u(t,x)\right )u(t,x)-ru(t,x)
\end{align*}
with
\begin{align*}
Q(u(t,x))&=\mathbbm{1}_{(-\infty ,v^{h})}(u(t,x))\gamma ^{h}+
\mathbbm{1}_{\left [ v^{l},\infty \right ) }(u(t,x))\gamma ^{l}
\\
&+\mathbbm{1}_{\left [ v^{h},v^{l} \right ) } (u(t,x))\left [
\frac{(\gamma ^{h}-\gamma ^{l})}{(v^{h}-v^{l})}\left (u(t,x)-v^{h}
\right )+\gamma ^{h}\right ].
\end{align*}
The terminal condition is $u(T,x)=g(x)$, where 
$g(x)=\min \{x_{1},\dots , x_{d}$\}.

\subsection{Hamilton-Jacobi-Bellman (HJB) equation}

Another type of nonlinear parabolic PDE is the HJB equation with its variables defined on Eq. (\ref{pde1}). In dynamic
programming, particularly in game theory with multiple players, each player
needs to solve the HJB equation to determine
an optimal strategy. Consider the stochastic dynamical model 
\begin{align}
\label{HJB10}
dX_{t}=2\sqrt{\lambda}\mathcal{M}_{t}dt+\sqrt{2} dW_{t},
\end{align}
where $X_{t}$ is a state process, $\{W_{t}\}$ is a Wiener process,
$\{\mathcal{M}_{t}\}$ is a control process, and $\lambda $ is a strength
of the control with $t\in [0,T]$ and $X_{0}=x\in \mathbbm{R}^{d}$. The objective of a control problem in a classical linear-quadratic
Gaussian (LQG) is to minimize a cost
functional
\begin{align*}
\displaystyle J\left (\{\mathcal{M}_{t}\}_{0\leq t\leq T}\right )=
\mathbbm{E}\left [\int _{0}^{T}||\mathcal{M}_{t}||^{2}dt+g(X_{T})
\right ].
\end{align*}

The aforementioned control problem can be written as an HJB equation for $d$ dimension 
\begin{align}
\frac{\partial u}{\partial t}(t,x)+\Delta _{x} u(t,x)=\lambda ||(
\nabla _{x}u)(t,x)||^{2}.
\label{HJB1}
\end{align}
Considering a terminal condition $u(T,x)=g(x)$ for
$x\in \mathbbm{R}^{d}$, an explicit solution of Eq. (\ref{HJB1}) can be
written as
\begin{align}
u(t,x)=-\frac{1}{\lambda} \ln{\left (\mathbbm{E}\left [\text{exp}(-
\lambda g(x+\sqrt{2}W_{T-t}))\right ]\right )}
\label{eq10}
\end{align}
with
\begin{align}
g(x)=\ln{\left (\frac{1+||x||^{2}}{2}\right )}.
\label{eq11}
\end{align}
For the case studied in this paper, we use the same number of assets in
\cite{han2018solving} and extend it to a higher dimension.
 
\subsection{The Deep-BSDE Method}

This section provides a brief overview of the deep-BSDE method, as outlined
in \cite{han2018solving}, to offer an insight about how the method works.
To numerically solve the BSDE in Eq. \eqref{3}, our primary objective
is to determine $u(0,X_{0})$, such that the value of
$u(t_{N},X_{t_{N}})$ calculated using Eq. \eqref{5} closely approximates
the terminal condition $g(X_{t_{N}})$. Additionally, calculation of the
value of $u(t_{N}, X_{t_{N}})$ depends on the unknown values of
$\nabla u(t_{n}, X_{t_{n}})$ for $n = 0, \ldots , N-1$. To estimate the
value of $u(0,X_{0})$, we adopt the deep-BSDE approach in~\cite{han2018solving}, which has been shown to solve BSDE better than
other traditional numerical methods. In general, we will update weight
parameters inside the deep neural networks as:
\begin{equation*}
{\boldsymbol{\theta}} = \left \{\theta _{u_{0}}, \theta _{\nabla u_{0}},
\theta _{1}, \ldots , \theta _{N-1} \right \},
\end{equation*}
where $\theta _{u_{0}}\approx u(0, X_{0})$,
$\theta _{\nabla u_{0}}\approx \nabla u(0,X_{0})$.

At the top level of deep neural network, the initial guess of
$\theta _{u_{0}}$ is selected randomly from a certain interval
$[a, b]$ which is set to be not far from the actual value
$u(0,X_{0})$. Based on Eqs. (\ref{4})-(\ref{6}), the deep neural
network estimates $u(t_{N}, X_{t_{N}})$ by optimizing the weight parameters
$\boldsymbol{\theta}$ in terms of a loss function $\ell $. In that case,
we minimize the loss function $\ell $, which is defined by
\begin{equation}
\label{lossfc}
\ell (\theta , \{ X_{t_{n}} \}, \{ \Delta W_{t_{n}} \})=\mathbb{E}\left ( ||g(X_{t_{N}})-u(t_{N},X_{t_{N}})||^{2}
\right ).
\end{equation}
The above loss function represents the difference of
$u(t_{N},X_{t_{N}})$ with the terminal condition $g(X_{t_{N}})$ of the
PDE. Then the deep-BSDE method optimizes all weights in the network using
the stochastic gradient descent-type (SGD) algorithm simultaneously. Note
that, the one we need the most is the solution of Eq. \eqref{pde1} or \eqref{3} which is represented by $\theta _{u_{0}}$ as an approximation
of $u(0,X_{0})$.

\color{black}

\section{Deep Genetic Algorithm (Deep-GA)}
\label{sec3}

In this research, we propose to combine GA with deep
learning techniques, resulting in a novel approach that we refer to as the deep
genetic algorithm (deep-GA). We embed the GA as an
optimization procedure into the well-cited deep-BSDE solver with the expectation
of enhancing its efficiency. Unlike other methods, such as the shooting
method \cite{burden2015numerical}, the GA is capable of handling the weights
of deep-BSDE and the stochastic differential equations. 
The embedding accelerates the convergence for calculating
$u_{0}$ while simultaneously widening the interval for determining the
initial guess $u_{0}$.

As mentioned earlier, the deep-BSDE method introduced in
\cite{han2018solving} requires an initial guess of $\theta _{u_{0}}$ as
the approximation of $u_{0}$, drawn from a certain interval $[a, b]$. This
$\theta _{u_{0}}$ represents the solution of the PDE, distinct from the
general weights of neural networks
$\left \{\theta _{\nabla u_{0}}, \theta _{1}, \ldots , \theta _{N-1}
\right \}$. Determining the appropriate value for the initial approximation,
$\theta _{u_{0}}$, proves to be a challenging task without precise knowledge
regarding the discussed problem's actual value. Our preliminary studies
(in Section~\ref{Prestudy}) found that the initial guess of
$\theta _{u_{0}}$ plays an important role in the solution process. When
$\theta _{u_{0}}$ is too far from the actual value, the method might need
prolonged iterations to converge or, in some cases, fail to converge altogether.
Consequently, finding a good initial guess is essential in the deep-BSDE
method.

\begin{figure*}
\includegraphics[width=\textwidth]{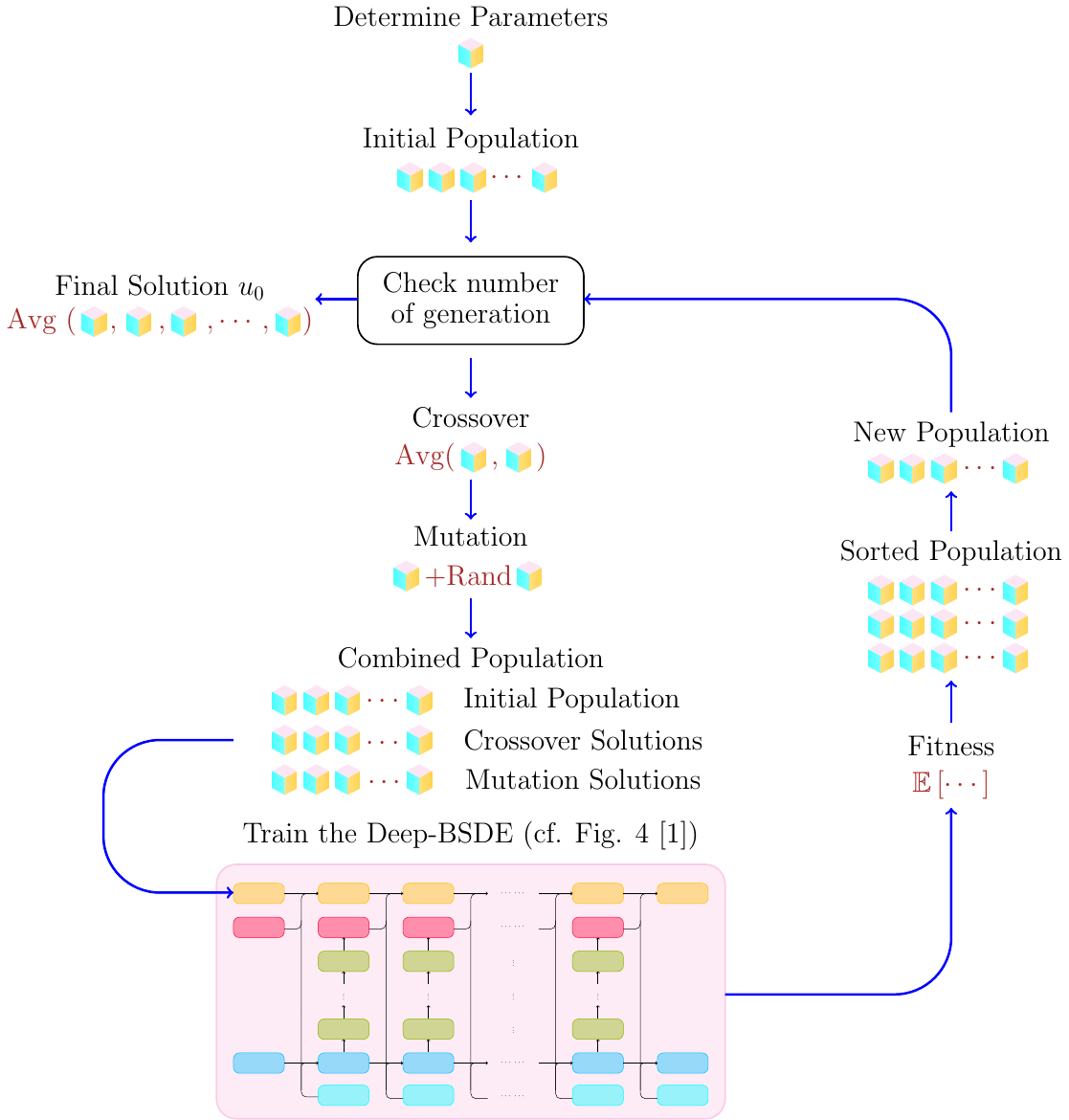}
\caption{Scheme of the deep-GA method.}
\label{figureDeepGA}
\end{figure*}

In our proposed method, we create a population of possible solutions to
approximate $u_{0}$ instead of using only one weight
$\theta _{u_{0}}$ in the original deep-BSDE method. The population is subsequently
optimized using the GA, while the remaining weights are updated using the
Adam optimizer akin to the process in the deep-BSDE method. A detailed
illustration of the deep-GA scheme is provided in Fig.~\ref{figureDeepGA}.

\subsection{Initialization and generating population}

The deep-GA begins by determining several important parameters, including
the number of chromosomes ($m$), the number of generations
($b$), the probability of crossover ($p_{c}$), the probability
of mutation ($p_{m}$), the minimum value of $u_{0}$ ($\min (u_{0})$), and
the maximum value of $u_{0}$ ($\max (u_{0})$). The parameters employed
in this study are $m=10$, $b=10$, $p_{c}=0.8$, $p_{m}=0.4$,
$\min (u_{0})=0$, $\max (u_{0})=100$ for the BS equation case,
and $\max (u_{0})=10$ for the HJB equation case. Notably, the interval
$[\min (u_{0}),\max (u_{0})]$ used for obtaining the initial guess of
$u_{0}$ is broader than the one used in the deep-BSDE method.

Next, a population $\{u_{0_{1}}, u_{0_{2}}, \dots , u_{0_{m}}\}$ consisting
of $m$ different possible solutions of $u_{0}$ is generated. The solutions
are selected from the interval $[\min (u_{0}),\max (u_{0})]$ using
\begin{align}
\label{npoint}
u_{0_{i}} = \min (u_{0}) + \frac{i}{m}(\max (u_{0})-\min (u_{0}))
\end{align}
for $i=1,2,\dots ,m$. Equation (\ref{npoint}) will give us a population consisting
of solutions which are evenly distributed over the interval. To ensure
stability in our results, three separate populations are generated in the
beginning of the deep-GA method. Genetic operators will evolve these populations
without any interactions among them. The final value of $u_{0}$ is taken
to be the average of all chromosomes from the three populations.

\subsection{Selection, crossover, and mutation operators}
To produce new approximations for $u_{0}$ from the population
$\{u_{0_{1}}, u_{0_{2}},\allowbreak \dots , u_{0_{m}}\}$, the GA uses
selection, crossover, and mutation operators. Each operator has a different
role in obtaining new approximations. In general, the selection is utilized
to pair solutions that will be used as inputs in the subsequent steps.
Crossover combines existing solutions to facilitate convergence within
a region, while mutation is employed to escape local optima
\cite{yang2020nature}.

The deep-GA method uses a random selection by choosing a pair of solutions
from the initial population. Let $(v_{1},v_{2})$ represent the selected
pair. With probability $p_{c}$, a new solution $v$ will be produced by
crossover which is calculated as
\begin{align}
\label{cross}
v = \frac{1}{2}(v_{1}+v_{2}).
\end{align}

Next, the mutation is applied to each solution $u_{0_{i}}$ in the initial
population with probability $p_{m}$ and a new solution $w$ produced is
calculated as
\begin{align}
\label{mutation}
w = u_{0_{i}} + \varepsilon \left (
\frac{\max (u_{0})-\min (u_{0})}{100} \right ),
\end{align}
where $\varepsilon $ is randomly selected from the interval $[-1,1]$.

In addition to the aforementioned mutations, we also compute four new solutions
namely
\begin{equation}
y_{j} = \bar{u}_{0} + \alpha _{j} \left (
\frac{\max (u_{0})-\min (u_{0})}{100} \right ),
\label{eq16}
\end{equation}
where $j=1,\dots ,4$, $\alpha _{1}=1$, $\alpha _{2}=-1$,
$\alpha _{3}=2$, $\alpha _{4}=-2$, and $\bar{u}_{0}$ is the average of
all solutions $u_{0_{i}}$, $i=1,\dots ,m$, in the initial population.

Selection, crossover, and mutation processes are performed $m$ times, corresponding
to the number of solutions in the initial population. The crossover and
mutation are controlled by two parameters $p_{c}$ and $p_{m}$ and the expected
number of new solutions produced is $[(p_{c}+p_{m})m]+4$. Combining with
the initial population, we would have about $[(p_{c}+p_{m}+1)m]+4$ solutions
in the current population.

\subsection{Deep neural network and fitness function}

A fitness function is required to assess the fitness value of each solution
$u_{0_{i}}$, $i=1,\dots ,m$, to determine which solutions are superior
to others. To define the fitness function, we use a similar deep neural
network with the same weights
$\boldsymbol{\theta} = \left \{\theta _{u_{0}}, \theta _{\nabla u_{0}},
\theta _{1}, \dots , \theta _{N-1}\right \}$ as those utilized in the deep-BSDE
method. In our proposed method, we optimize $\theta _{u_{0}}$ using the
GA and optimize
$\left \{\theta _{\nabla u_{0}}, \theta _{1}, \dots , \theta _{N-1}
\right \}$ using an SGD-type algorithm. The optimization processes are
carried out alternately.

First, we set the weight $\theta_{u_{0}}$ with the average of all solutions
from the current three populations. The network is then trained for
$p$ times (iterations) in each generation using the Adam optimizer. Even
though at this step $\theta _{u_{0}}$ is also updated, the value will be
discarded and not used. In order to calculate the fitness value of
$u_{0_{i}}$, we assign and replace the weight $\theta _{u_{0}}$ in the
network with $u_{0_{i}}$. Then the fitness function of the deep-GA is defined
similarly to Eq. \eqref{lossfc}, but with a conditional weight as follows,
\begin{equation}
\label{fitnessfc}
f(u_{0_{i}}, \boldsymbol{\theta}, \{ X_{t_{n}} \}, \{ \Delta W_{t_{n}}
\}) = \mathbb{E}(||g(X_{t_{N}})-u(t_{N},X_{t_{N}})||^{2} \;\;|\;\theta _{u_{0}}=u_{0_{i}}).
\end{equation}
When calculating the fitness values of $u_{0_{i}}$, for
$i=1,\dots ,m$, the other weights
$\left \{\theta _{\nabla u_{0}}, \theta _{1}, \dots , \theta _{N-1}
\right \}$ remain unchanged.
In this work, the fitness function is the same as that from Han et al.\ \cite{han2018solving}, i.e., Eq. (\ref{lossfc}). The conditional weights in
the fitness function are selected to determine the value of
$\theta _{u_{0}}$ and obtain a precise fitness value for
$u_{0_{i}}$.

The deep-GA method only uses a single deep neural network and shares its
weights across all solutions $u_{0_{i}}$ in the population. These processes
are repeated in each generation. Since the deep-GA method consists of
$b$ generations, the network is trained for a total of $pg$ iterations. The
value of $p$ can be chosen from $[100,1000]$.

\subsection{New population and final solution}
 
After completing all the aforementioned stages, it is necessary to sort
the solutions in the population based on their fitness values from Eq. \eqref{fitnessfc}.
The top $m$ solutions will continue as the new population, while the remaining
solutions are eliminated. All the processes in the deep-GA method are repeated
for $b$ generations. Finally, the solution $u_{0}$ of the BSDEs, obtained
by the deep-GA method, is the average of all solutions from the three populations.

Pseudo-code and scheme of the deep-GA method are illustrated in Fig.~\ref{figureDeepGA}
and Algorithm~\ref{alg:two}. Note that for simplicity, the scheme in Fig.~\ref{figureDeepGA} is only using one population, instead of three.
The code used for calculating the fitness function in the deep-GA is sourced
from \cite{han2018solving,deepBSDE}.
Our code is based on the BSDE solver \cite{han2018solving} and is available at \cite{GABSDE}. We compiled the deep-BSDE and the deep-GA using Google Collaboratory.

\begin{algorithm}
\caption{Pseudo-code of the deep-GA method.}\label{alg:two}
\textbf{Input:} $m$, $b$, $p_{c}$, $p_{m}$, $\min (u_{0})$, $\max (u_{0})$,
$\alpha _{1}$, $\alpha _{2}$, $\alpha _{3}$, $\alpha _{4}$ \;
\textbf{Output:} $u_{0}$ \;

$P_{1} \gets $ GeneratePopulation($\max (u_{0}), \min (u_{0}), m$) \;
$P_{2} \gets $ GeneratePopulation($\max (u_{0}), \min (u_{0}), m$) \;
$P_{3} \gets $ GeneratePopulation($\max (u_{0}), \min (u_{0}), m$) \;

\For{$\textnormal{generation} = 1: b$}{
    \For{$l = 1: 3$}{
        $\bar{u}_{0} \gets $ Average($P_{l}$) \;
        \For{$i = 1: m$}{
            \If{$\textnormal{Random}([0,1]) < p_{c}$}{
                $v_{1} \gets $ RandomSelection($P_{l}$) \;
                $v_{2} \gets $ RandomSelection($P_{l}$) \;
                $v \gets \frac{1}{2}(v_{1} + v_{2})$ \;
                $P_{l} \gets $ Append($P_{l},v$) \;
            }
            \If{$\textnormal{Random}([0,1]) < p_{m}$}{
                $\varepsilon \gets \textnormal{Random}([-1,1])$ \;
                $w \gets u_{0_{i}} + \varepsilon \left ( \frac{\max (u_{0})-\min (u_{0})}{100} \right )$ \;
                $P_{l} \gets $ Append($P_{l},w$) \;
            }
        }
        \For{$j = 1: 4$}{
            $y_{j} \gets \bar{u}_{0} + \alpha _{j} \left ( \frac{\max (u_{0})-\min (u_{0})}{100} \right )$ \;
            $P_{l} \gets $ Append($P_{l},y_{j}$) \;
        }
    }
    $\theta _{u_{0}} \gets $ Average($P_{1},P_{2},P_{3}$) \;
    \For{$i = 1: p$}{
        $\{ X_{t_{n}} \}, \{ \Delta W_{t_{n}} \} \gets $ GenerateSample($\mu ,\sigma ,\Delta t$) \;
        $\theta \gets $ Train($\theta , \{ X_{t_{n}} \}, \{ \Delta W_{t_{n}} \}$) \;
    }
    $\{ X_{t_{n}} \}, \{ \Delta W_{t_{n}} \} \gets $ GenerateSample($\mu ,\sigma ,\Delta t$) \;
    \For{$l = 1: 3$}{
        \For{$i = 1: \textnormal{Size}(P_{l})$}{
            $\theta _{u_{0}} \gets u_{0_{i}}$ \;
            fitness $\gets f(u_{0_{i}},\theta , \{ X_{t_{n}} \}, \{ \Delta W_{t_{n}} \})$ \;
            $F_{l} \gets $ Append($F_{l}, \textnormal{fitness}$) \;
        }
    }
    $P_{l}, F_{l} \gets $ Sort($P_{l}, F_{l}$)\;
    $P_{l} \gets $ Eliminate($P_{l}, m$)\;
}
$u_{0} \gets $ $\frac{1}{3}$(Average($P_{1}$) $+$ Average($P_{2}$) $+$ Average($P_{3}$)) \;
\end{algorithm}

\section{Results and Discussions}
\label{sec4}

\subsection{Preliminary study}\label{Prestudy}

\begin{figure}
\centering
\includegraphics[width=0.8\textwidth]{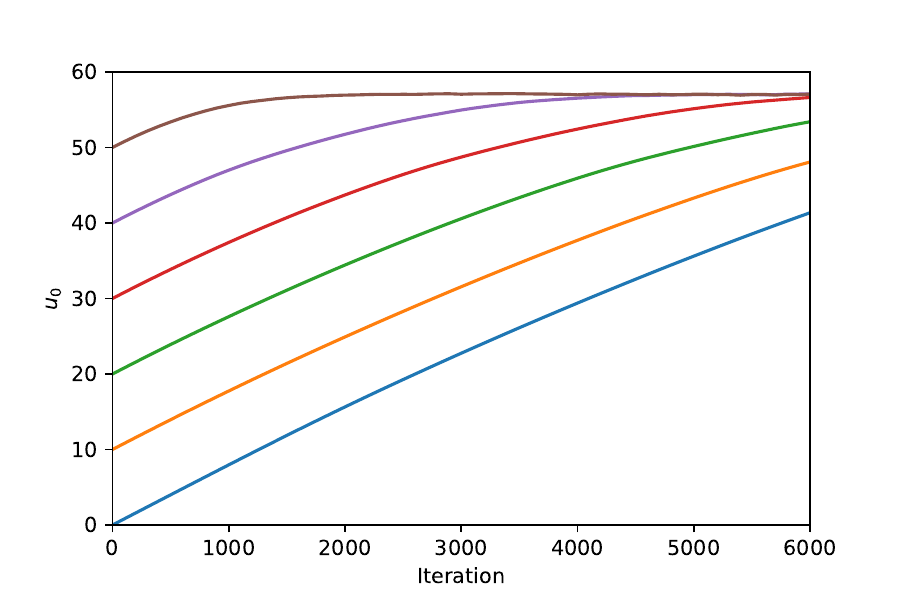}
\caption{Plot of $u_{0}$ obtained from the deep-BSDE method in 6000 iterations
for six different initial guesses of $u_{0}$.}
\label{FigGA3}
\end{figure}

First, we demonstrate the search for solutions of the BS equation
using the deep-BSDE method with various initial guesses of
$\theta _{u_{0}}$. All settings used in the simulation follow
\cite{han2018solving} and the network is trained for 6000 iterations. Figure \ref{FigGA3} shows that when the initial guess of $\theta _{u_{0}}$ is
significantly far from the actual value, the convergence of
$u_{0}$ is slow. In fact, convergence cannot be achieved within 6000 iterations
for such initial guesses. Thus, the selection of the initial guess plays
a crucial role in the solution process.

\begin{table*}
\centering\small
\caption{Loss values of the deep-BSDE method for different initial guesses of
$\theta _{u_{0}}$ from five independent runs without any weight training - the
nonlinear - BS equation with default risk.}
\label{tab01}
\begin{tabular}{crrrrrrr}
\hline
\multirow{2}{*}{$\theta_{u_0}$} & \multicolumn{5}{c}{Loss} & \multirow{2}{*}{Average} & \multirow{2}{*}{Std} \\ \cline{2-6}
& \multicolumn{1}{c}{1} & \multicolumn{1}{c}{2} & \multicolumn{1}{c}{3} & \multicolumn{1}{c}{4} & \multicolumn{1}{c}{5} & & \\ \hline
0 & 3682.33 & 3581.65 & 3678.51 & 3677.77 & 3723.98 & 3668.85 & 52.43 \\
10 & 2491.00 & 2488.57 & 2488.68 & 2466.36 & 2503.04 & 2487.53 & 13.26 \\
20 & 1525.25 & 1564.84 & 1560.57 & 1558.50 & 1526.19 & 1547.07 & 19.63 \\
30 & 786.31 & 796.03 & 787.36 & 815.15 & 815.17 & 800.00 & 14.34 \\
40 & 312.10 & 330.30 & 324.18 & 307.54 & 323.30 & 319.48 & 9.36 \\
50 & 72.55 & 67.43 & 72.57 & 67.86 & 64.30 & 68.94 & 3.58 \\
60 & 34.30 & 40.71 & 27.08 & 32.05 & 28.79 & 32.59 & 5.34 \\
70 & 146.43 & 147.19 & 155.63 & 148.78 & 173.33 & 154.27 & 11.25 \\
80 & 473.34 & 484.84 & 513.87 & 487.00 & 511.87 & 494.18 & 17.84 \\
90 & 1013.13 & 1012.78 & 1013.73 & 1079.93 & 1068.09 & 1037.53 & 33.56 \\
100 & 1772.72 & 1812.55 & 1782.69 & 1824.90 & 1798.66 & 1798.31 & 21.26 \\ \hline
\end{tabular}
\end{table*}

\begin{table}
\centering\small
\caption{Loss values of the deep-BSDE method for different initial guesses of
$\theta _{u_{0}}$ from five independent runs without any weight training - the
HJB equation.}
\label{tab01a}
\begin{tabular}{crrrrrrr}
\hline
\multirow{2}{*}{$\theta_{u_0}$} & \multicolumn{5}{c}{Loss} & \multirow{2}{*}{Average} & \multirow{2}{*}{Std} \\ \cline{2-6}
& \multicolumn{1}{c}{1} & \multicolumn{1}{c}{2} & \multicolumn{1}{c}{3} & \multicolumn{1}{c}{4} & \multicolumn{1}{c}{5} & & \\ \hline
0 & 19.821 & 19.512 & 19.565 & 19.485 & 19.992 & 19.675 & 0.198 \\
1 & 11.725 & 11.821 & 11.941 & 11.536 & 11.732 & 11.751 & 0.133 \\
2 & 5.801 & 6.149 & 5.965 & 6.037 & 6.021 & 5.995 & 0.114\\
3 & 2.086 & 2.034 & 2.096 & 2.168 & 2.087 & 2.094 & 0.043 \\
4 & 0.232 & 0.208 & 0.225 & 0.231 & 0.262 & 0.232 & 0.017 \\
5 & 0.414 & 0.360 & 0.349 & 0.420 & 0.348 & 0.378 & 0.032 \\
6 & 2.492 & 2.511 & 2.483 & 2.535 & 2.564 & 2.517 & 0.030 \\
7 & 6.575 & 6.526 & 6.588 & 6.721 & 6.674 & 6.617 & 0.071 \\
8 & 12.746 & 13.037 & 12.840 & 12.500 & 12.770 & 12.779 & 0.173 \\
9 & 20.578 & 21.118 & 20.889 & 20.684 & 21.128 & 20.879 & 0.223 \\
10 & 31.104 & 30.949 & 31.036 & 31.141 & 31.015 & 31.049 & 0.068 \\ \hline
\end{tabular}
\end{table}

Secondly, we illustrate the significance of selecting initial guesses for
$\theta _{u_{0}}$ during the first iteration of deep-BSDE. Each initial
guess is compiled in five independent runs to demonstrate that the random
weights of deep-BSDE in each run have minimal impact on the loss values,
while the choice of $\theta _{u_{0}}$ does. The results are presented in
Table~\ref{tab01} and Table~\ref{tab01a}. The loss values
of all independent runs have relatively small changes compared to the averages
for each $u_{0}$ value. As indicated by the standard deviations, when the
initial guess $\theta _{u_{0}}$ is closer to the actual values, the loss values are smaller. This indicates that to enhance efficiency, it is important to select
$\theta _{u_{0}}$ with greater accuracy.

\subsection{Non-linear Black-Scholes equation with default risk}

The fair price of the European multi asset contingent claim, as governed
by the nonlinear BS equation with default risk in Eq. (\ref{3}), has
been simulated using both our proposed deep-GA and the deep-BSDE methods. For the simulations, we use
parameters in \cite{han2018solving}, i.e., $r=0.02$,
$\delta =2/3$, $\gamma _{h}=0.2$, $\gamma _{l}=0.02$,
$\hat{\mu}=0.02$, $v_{h}=50$, and $v_{l}=70$. We assume that 
$d=100$ assets have uniform price $x=\{100,100,...,100\}$ with uniform
volatility values $\hat{\sigma}$ for each asset. It is important to note
that the volatility $\sigma $ is a $d\times d$ matrix for each value (see Eq. (\ref{pde1})) but here we use $\hat{\sigma}$ for simulation based on
Eq. (\ref{Eq.2}). The connection between $\sigma $ and
$\hat{\sigma}$ has been explained in Section~\ref{NLBS}. In addition, we range
the values of volatility used in the simulation
$\hat{\sigma}\in \left \{0.1,0.5\right \}$ to show the effect of volatility
to the value of the European contingent claim. The parameters of GA are
also determined following \cite{han2018solving}: learning rate = 0.008,
batch size = 64, and validation size = 256. For the deep-BSDE method, we
choose a single value from an interval $[40,50]$ as the initial guess value
of $\theta _{u_{0}}$ and the model is trained for 10000 times. For the
deep-GA method, we use $g=15$, $p=100$, and hence the model is trained
1500 times in total.

We compare the results obtained by the deep-GA and the deep-BSDE methods for various values of $\hat{\sigma}$ in {Table~\ref{tab03}}. Following \cite{weinan2017deep,weinan2019multilevel}, we will take the values of $u_{0}$ obtained using Picard method as ``the actual values''. The
results from the deep-GA method exhibit strong agreement with the ``actual
values'', and the absolute percentage error
is comparable to that of the deep-BSDE.

\begin{table*}
	\centering\small
\caption{Results obtained by the deep-BSDE and the deep-GA for different values
of $\hat{\sigma}$.}
\label{tab03}
\begin{tabular}{cccccccc}
		\hline
		\multirow{2}{*}{$\hat{\sigma}$} & \multicolumn{3}{c}{$u_0$} & \multicolumn{2}{c}{Abs. Percentage Error} & \multicolumn{2}{c}{Time}\\ \cline{2-4}\cline{5-6}\cline{7-8}
		& {\footnotesize Picard} & {\footnotesize Deep-BSDE} & {\footnotesize Deep-GA} & {\footnotesize Deep-BSDE} & {\footnotesize Deep-GA} & {\footnotesize Deep-BSDE} & {\footnotesize Deep-GA} \\ \hline
    0.1 & 77.00 & 76.88 & 76.95 & 0.156 & 0.065 & 1386 & 759 \\
    0.2 & 57.32 & 56.99 & 57.07 & 0.576 & 0.436 & 1351 & 736 \\
    0.3 & 42.50 & 42.24 & 42.05 & 0.612 & 1.059 & 1350 & 749 \\
    0.4 & 32.12 & 31.63 & 31.94 & 1.526 & 0.560 & 1348 & 775 \\
    0.5 & 24.09 & 23.50 & 23.38 & 2.449 & 2.947 & 1436 & 775 \\ \hline
    \end{tabular}
\end{table*}

Various values of volatility $\hat{\sigma}$ represent changes in assets'
rate of returns, leading to increased complexity in the pricing process.
The price of contingent claim $u(x,t)$ is determined by solving Eq. (\ref{3}).
As $\hat{\sigma}$ increases, the task determining the initial guess of $\theta _{u_{0}}$ becomes increasingly complicated, resulting in longer computational times. However, our deep-GA shows a higher efficiency than that of the deep-BSDE method as the volatility $\hat{\sigma}$ increases.

\begin{table*}
	\centering\small
\caption{Comparison of the deep-BSDE and the deep-GA for different numbers of
dimension - the nonlinear BS equation with default risk.}
\label{tab04}
\begin{tabular}{cccccccc}
    \hline
		\multirow{2}{*}{$d$} & \multicolumn{3}{c}{$u_0$} & \multicolumn{2}{c}{Abs. Percentage Error} & \multicolumn{2}{c}{Time }\\ \cline{2-4}\cline{5-6}\cline{7-8}
		& {\footnotesize Picard} & {\footnotesize Deep-BSDE} & {\footnotesize Deep-GA} & {\footnotesize Deep-BSDE} & {\footnotesize Deep-GA} & {\footnotesize Deep-BSDE} & {\footnotesize Deep-GA} \\ \hline
    100 & 77.00 & 76.88 & 76.95 & 0.156 & 0.065 & 1386 & 759 \\
		200 & 75.16 & 75.08 & 75.24 & 0.106 & 0.093 & 2315 & 1064 \\
		300 & 74.18 & 74.12 & 74.06 & 0.081 & 0.162 & 3340 & 1295 \\
		400 & 73.53 & 73.48 & 73.41 & 0.068 & 0.163 & 4438 & 1444 \\
		500 & 72.99 & 72.97 & 72.98 & 0.027 & 0.014 & 5571 & 1632 \\\hline
    \end{tabular}
\end{table*}

To show the high-dimensional effect to the solution process of the BS
equation, we present some numerical results in {Table~\ref{tab04}}. As the
dimension increases, the computational time for both methods increases.
Nevertheless, the deep-GA has a half rate of computational time increment
than that of the deep-BSDE in average. Therefore, the deep-GA is more efficient
in solving the PDE than the deep-BSDE.%

\begin{figure}
\centering
\includegraphics[width=0.8\textwidth]{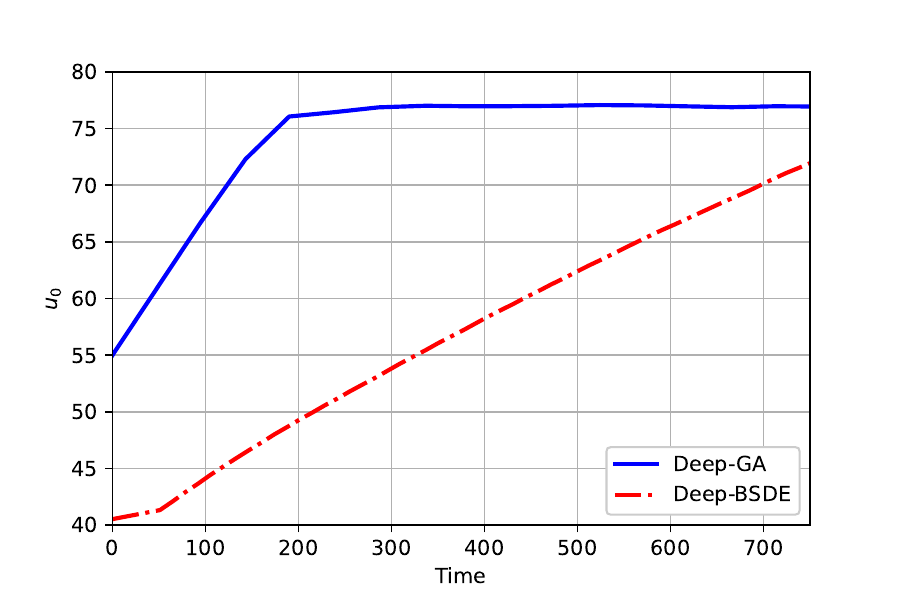}
\caption{$u_{0}$ of the BS equation against time for
$\hat{\sigma}=0.1$ in the first 750 seconds.}
\label{figure11a}
\end{figure}

\begin{figure}
\centering
\includegraphics[width=0.8\textwidth]{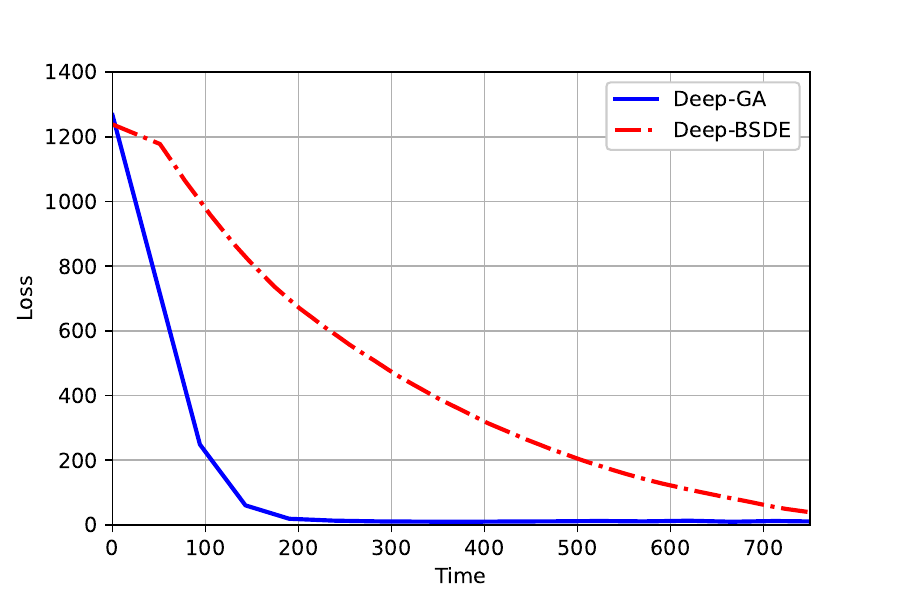}
\caption{Loss in the BS equation against time for $\hat{\sigma}=0.1$
in the first 750 seconds.}
\label{figure11b}
\end{figure}

Furthermore, we present in Figs.~\ref{figure11a} and \ref{figure11b} the convergence of the estimated solutions $u_{0}$ and the losses for $\hat{\sigma} =0.1$ in the first 750 seconds of both methods, respectively. In the case of the deep-GA method, the solutions $u_{0}$ and the loss values plotted in the figures represent
the average of all solutions and the averages of all fitness values accross
the three populations. One can observe the fast convergent of the solution computed using our proposed method. Although losses from both methods
started with nearly identical initial values, the loss of the deep-GA method decreases rapidly in the first 200 seconds.

\subsection{Hamilton-Jacobi-Bellman (HJB) equation}

In this subsection, we simulate an HJB equation
for the minimum cost of an optimal strategy in  investments. For simulations,
we use $\lambda \in \left \{1,10,20,30,40,50\right \}$, learning rate =
0.01, batch size = 64, and validation size = 256. For the deep-BSDE method,
we choose a random value from the interval $[7,8]$ as an initial guess
of $\theta _{u_{0}}$ and then the model is trained for 40000 times, while
for the deep-GA method, we use $b=20$, $p=1000$, resulting in a total of
20,000 training iterations.

\begin{table*}
	\centering\small
\caption{Results obtained by the deep-BSDE and the deep-GA for different values
of $\lambda $.}
\label{tab002}
\begin{tabular}{cccccccc}
		\hline
		\multirow{2}{*}{$\lambda$} & \multicolumn{3}{c}{$u_0$} & \multicolumn{2}{c}{Abs. Percentage Error} & \multicolumn{2}{c}{Time}\\ \cline{2-4}\cline{5-6}\cline{7-8}
		& {\footnotesize Monte Carlo} & {\footnotesize Deep-BSDE} & {\footnotesize Deep-GA} & {\footnotesize Deep-BSDE} & {\footnotesize Deep-GA} & {\footnotesize Deep-BSDE} & {\footnotesize Deep-GA} \\ \hline
    1  & 4.590 & 4.606 & 4.596 & 0.349 & 0.131 & 2531 & 1530 \\
    10 & 4.493 & 4.503 & 4.508 & 0.223 & 0.334 & 2680 & 1512 \\
    20 & 4.369 & 4.415 & 4.416 & 1.053 & 1.076 & 2658 & 1495 \\
    30 & 4.247 & 4.370 & 4.350 & 2.896 & 2.425 & 2691 & 1483 \\
    40 & 4.158 & 4.295 & 4.281 & 3.295 & 2.958 & 2645 & 1506 \\
    50 & 4.096 & 4.241 & 4.212 & 3.540 & 2.832 & 2706 & 1504 \\ \hline
    \end{tabular}
\end{table*}

The results obtained from both the deep-BSDE and deep-GA methods for various
values of the control strength $\lambda $ are presented in {Table~\ref{tab002}}. We compare these results with those obtained using the Monte
Carlo method, as reported in \cite{han2018solving}. Assuming Monte Carlo
method provides ``actual results'', the solution obtained by the deep-GA
is fairly comparable to that obtained by the deep-BSDE method but with
significantly less computational time. The use of the GA in the deep-GA method for optimizing the search of the initial guess
value noticeably enhances the efficiency compared to the former method
described in \cite{han2018solving}.

\begin{table*}
	\centering\small
\caption{Comparison of the deep-BSDE and the deep-GA for different numbers of
dimension - the HJB equation.}\label{tab02}
\begin{tabular}{cccccccc}
    \hline
		\multirow{2}{*}{$d$} & \multicolumn{3}{c}{$u_0$} & \multicolumn{2}{c}{Abs. Percentage Error} & \multicolumn{2}{c}{Time}\\ \cline{2-4}\cline{5-6}\cline{7-8}
		& {\footnotesize Monte Carlo} & {\footnotesize Deep-BSDE} & {\footnotesize Deep-GA} & {\footnotesize Deep-BSDE} & {\footnotesize Deep-GA} & {\footnotesize Deep-BSDE} & {\footnotesize Deep-GA} \\ \hline
    100 & 4.590 & 4.606 & 4.596 & 0.349 & 0.131 & 2531 & 1530 \\
		200 & 5.291 & 5.296 & 5.296 & 0.095 & 0.095 & 4225 & 1740 \\
		300 & 5.699 & 5.706 & 5.697 & 0.123 & 0.035 & 5776 & 2406 \\
		400 & 5.988 & 5.983 & 5.985 & 0.084 & 0.050 & 7695 & 2523 \\
		500 & 6.212 & 6.220 & 6.211 & 0.129 & 0.016 & 9859 & 2943 \\
    \hline
    \end{tabular}
\vspace{6pt}
\end{table*}

For various numbers of dimension $d$, the computational times of both deep-GA
and deep-BSDE to solve the HJB equation increase as the dimension is higher.
The rate of computational time increase of the deep-GA is also approximately
half of the deep-BSDE in average. The results for the HJB equation are
consistent with those for the BS equation (see Eqs.(\ref{Eq.2})
and (\ref{HJB1})). Moreover, the percentage error of the deep-GA also follows the same trend on average with that of the deep-BSDE (see Table~\ref{tab02}) .

\begin{figure}
\centering
\includegraphics[width=0.8\textwidth]{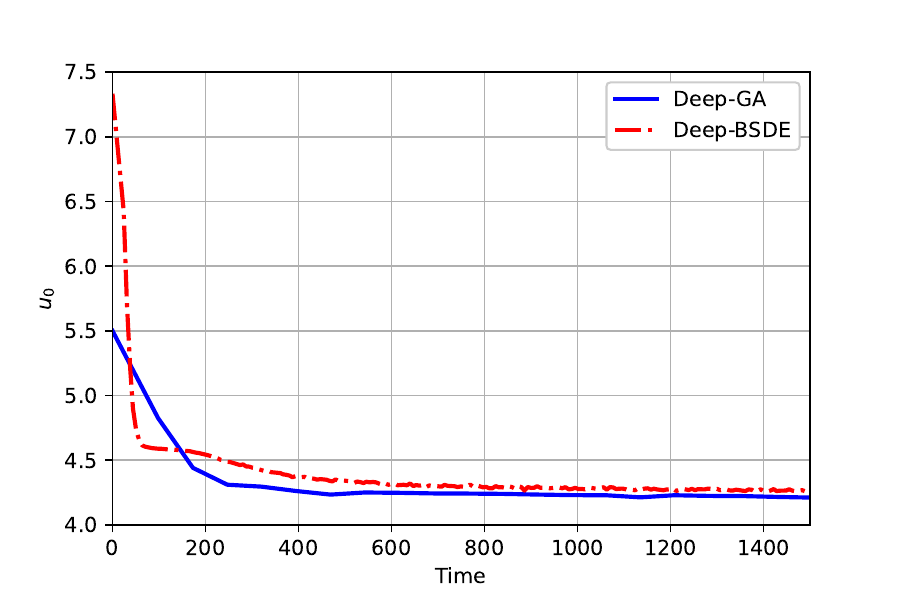}
\caption{$u_{0}$ of the HJB equation against time for $\lambda =50$ in the first
1500 seconds.}
\label{figure12a}
\end{figure}

\begin{figure}
\centering
\includegraphics[width=0.8\textwidth]{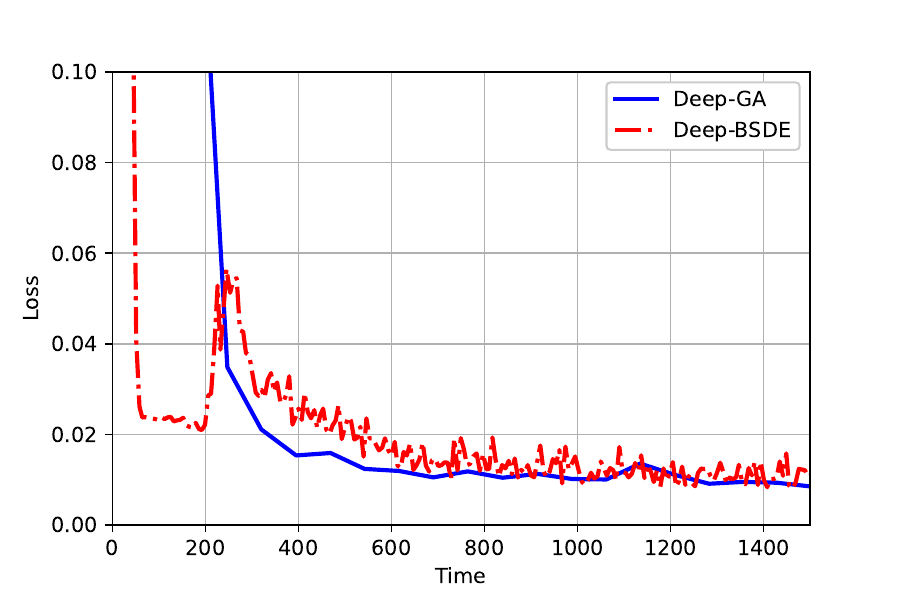}
\caption{Loss of the HJB equation against time for $\lambda =50$ in the first
1500 seconds.}
\label{figure12b}
\end{figure}

In Figs.~\ref{figure12a} and \ref{figure12b}, we show the trajectories
of the estimated solutions $u_{0}$ and the associated losses for
$\lambda =50$ in the first 1500 seconds from both methods. Notably, the
solutions of the deep-BSDE method decrease faster only in the first 150
seconds. Beyond that, the deep-GA method demonstrates superiority and
rapid convergence. This pattern is also reflected in the movement of losses
for both methods. Despite using the same learning rate, the loss of the
deep-BSDE method displays noticeable fluctuations while the one of the
deep-GA methods remains stable. Full detailed results of different
parameters from both methods for the BS and HJB cases are provided
at \cite{GABSDE}.

\section{Conclusion}
\label{sec5}

In this research, we introduced a novel method for solving high-dimensional
and nonlinear PDEs, namely the deep-GA
method. We applied this method to estimate the solution $u_{0}$ of BSDEs associated with two nonlinear PDEs: the BS equation with default risk and the HJB equation, considering various parameters. The results obtained using
this method for both equations demonstrate good agreement with the actual
values obtained from Picard and Monte Carlo methods. Furthermore, the deep-GA
method exhibits significantly efficient computational times compared
to the deep-BSDE method.

We have demonstrated that giving careful consideration to the choice of the
initial guess for $\theta _{u_{0}}$ rather than solely focusing on
updating weights of the hidden layers is critical. This approach provides
a broad and flexible way to determine the initial guess interval,
resulting in a fast convergence rate compared to the deep-BSDE method.
Furthermore, our proposed method can offer effective suggestions even when
the ideal interval for the initial guess of $\theta _{u_{0}}$ is unknown.
We believe that this method can be applied to solve other problems where
one or more weights play a more significant role than others.
For future work, a jump-diffussion case might be considered to represent the existence of
small and continuous price changes and large and infrequent ones represented
by Poisson process, simultaneously.

\section{Acknowledgement}
This work was funded by the Ministry of Education, Culture, Research, and Technology of the Republic Indonesia through the World Class Professor Program year 2021. H.S.\ also acknowledges support from Khalifa University through a Faculty Start-Up Grant (No.\ 8474000351/FSU-2021-011) and a Competitive Internal Research Awards Grant (No.\ 8474000413/CIRA-2021-065).

\bibliographystyle{elsarticle-num} 
\bibliography{main_document}

\end{document}